\documentclass[12pt,a4paper]{article}
\usepackage{amsmath,amssymb,amsfonts}
\usepackage[latin1]{inputenc}
\def\Bbb{\mathbb}

\def\BQ{``}
\def\EQ{" }

\title{\bf Power residues, digit expansions and relative class numbers}

\author{Kurt Girstmair}

\date{}

\makeatletter
\let\@@maketitle=\maketitle
\def\maketitle{\def\thispagestyle##1{\relax}\@@maketitle}
\makeatother
\newtheorem{theorem}{Theorem}

\newtheorem{lemma}{Lemma}
\newtheorem{corollary}{Corollary}

\def\BE{\begin{equation}}
\def\EE{\end{equation}}
\def\BD{\begin{displaymath}}
\def\ED{\end{displaymath}}
\def\BA{\begin{array}}
\def\EA{\end{array}}
\def\BEA{\begin{eqnarray*}}
\def\EEA{\end{eqnarray*}}
\def\BI{\bibitem}

\def\Z{\Bbb Z}
\def\Q{\Bbb Q}
\def\R{\Bbb R}
\def\C{\Bbb C}

\def\phi{\varphi}

\def\MB{\mbox}
\def\LD{\ldots}
\def\OV{\overline}
\def\SP#1{\langle #1 \rangle}
\def\WH{\widehat}
\def\WT{\widetilde}

\def\SIGN{\MB{Sign}}

\def\sminus{\smallsetminus}

\def\NDIV{\, \nmid \,}

\def\MN{\medskip\noindent}
\def\STOP{\hfill$\Box$}

\def\B3{B^{(3)}}
\def\ct3h{\WH{\rm{ct}}^{(3)}}

\def\LS#1#2{ \left( \frac{#1}{#2} \right) }
\def\SB{S^{(b)}}
\def\SZ{S^{(2)}}
\def\SPP{S^{(p+1)}}
\def\EB{E^{(b)}}
\def\TB{T^{(b)}}

\def\BCH{B_{\chi}}
\def\BCHB{B_{\OV{\chi}}}
\def\OVV#1{\OV{\OV{#1}}}

\def\SIGN{\MB{sign}}

\begin{document}
\maketitle

\begin{abstract}
\noindent
This is a survey of a connection between the distribution of certain power residues modulo $p$, $p$ a prime, and relative class numbers. The focus lies on quadratic residues and
sixth power residues. Dirichlet's class number formula yields a number of results about the distribution of quadratic residues, for instance, the well-known fact that the interval
$[0,p/2]$ contains more quadratic residues than nonresidues. This class number formula is also responsible for some properties of the digit expansions of numbers $m/p$, $p\NDIV m$.
In a certain sense the results based on Dirichlet's formula can be extended to sixth power residues, where geometry plays an important role.
\end{abstract}

\section*{1. Introduction}

Let $p$ be a prime $\equiv 3$ mod $4$, $p>3$. Let $C$ denote the set of quadratic residues mod $p$ in $\{1,2,\LD,p-1\}$ and $N$ the set of quadratic nonresidues in this set.
We have $|C|=|N|=(p-1)/2$. The following result is due to Dirichlet:
\BE
\label{1.1}
  |C\cap[0,p/2]|-|N\cap[0,p/2]|>0;
\EE
see \cite[p. 373]{BoSa}
Hence there are more quadratic residues in the interval $[0,p/2]$ than quadratic nonresidues.
More precisely,
 \BE
\label{1.2}
   |C\cap[0,p/2]|-|N\cap[0,p/2]| =\begin{cases}
                                    h, \hbox{ if $p\equiv 7$ mod 8};\\
                                  3h, \hbox{  if $p\equiv 3$ mod 8},
                                 \end{cases}
\EE
where $h$ is the class number of the imaginary quadratic number field $\Q(\sqrt{-p})$. So this class number determines not only sign of the difference $|C\cap[0,p/2]|-|N\cap[0,p/2]|$,
but also its magnitude.

From (\ref{1.2}) one obtains the following result, which is shown in \cite{Be}.
\BE
\label{1.3}
   |C\cap[0,p/6]|-|N\cap[0,p/6]|=\begin{cases}
                                 -h, \hbox{ if $\LS 2p=\LS 3p=-1$},\\
                                  h, \hbox{ otherwise}.
                              \end{cases}
\EE
Here $\LS{-}p$ is the Legendre symbol.
Accordingly, the class numbers determines the absolute value of this difference, whereas its sign depends on $\{2,3\}\subseteq N$.

Moreover, the said class number formula implies the following result of \cite{Gi1}.
For $p\equiv 3$ mod 4, $p>3$,  and $b\ge 2$ let
\BE
\label{1.3.1}
  \frac 1p=\sum_{k=1}^{\infty}a_kb^{-k}
\EE
be the {\em digit expansion} of $1/p$ with respect to the basis $b$, i.e., the digits $a_k$ take only values in $\{0,1,\LD,b-1\}$.
If, in addition, $b$ is a primitive root mod $p$, then $(a_1,a_2,\LD,a_{p-1})$ is a period of this expansion. Then we have
\BE
\label{1.4}
  (a_2+a_4+\LD+a_{p-1})-(a_1+a_3+\LD+a_{p-2})=(b+1)h,
\EE
where $h$ is as above. This means that the
difference on the left-hand side of (\ref{1.4}) is always positive. Therefore, the
sum of the digits with even indices in the period exceeds the sum of the digits with odd indices.
For example, 10 is a primitive root mod 7 and $1/7=0.\OV{142857}$,  where the bar marks the period. Then  $(4+8+7)-(1+2+5)=11=(10+1)h$, since $h=1$ for $p=7$.

A common viewpoint for the results (\ref{1.2}), (\ref{1.3}) and (\ref{1.4}) was developed in \cite{Gi2}. Furthermore, the respective statements concerning the signs of the
left-hand sides of these identities were generalized in a certain sense, in particular, to sixth power residues.
However, the generality of this approach has the effect that the most interesting special cases get lost in some sense,
in particular, since their description is scattered over several sections of a long paper.

In contrast to this situation, the present article is focused on the most relevant special cases, which become clearer in this way, as we hope.
But it is not merely an epitome of the said paper. It also contains new developments
like Theorems \ref{s2} and \ref{s2.1} on digit expansions, which considerably generalize known results.
In the case of sixth power residues, the two diagrams are also novel and allow new insights.
Moreover, we illustrate the results by many examples.
Finally, we point out what is possible for tenth power residues and hint on some more recent developments in the literature.

\section*{2. The $b$-deviation vector} 

In what follows let $p$ be a prime and $q$ an even divisor of $p-1$. The cyclic group $G=(\Z/p\Z)^{\times}$ consists of the residue classes $\OV 1, \OV 2, \LD, \OV{p-1}$.
It has exactly one subgroup of index $q$, namely $H=\{\OV k^q; k=1,\LD,p-1\}$. We require $\OV{-1}\not\in H$, which means $p\equiv q+1$ mod $2q$.
The factor group $G/H$ consists of the classes $\OVV k$, $k=1,2,\LD,p-1$. We identify each class
$C\in G/H$ with $\{k; 1\le k\le p-1, \OVV k =C\}$. If $q=2$, then $G/H$ consists of the classes $H$ and $\{ 1,2,\LD,p-1\}\sminus H$,
which comprise the quadratic residues and nonresidues mod $p$, respectively. As a rule, we describe the classes $C\in G/H$ as follows.
Let $g$ be a primitive root $p$. Then
\BD
G/H=\{\OVV g^{\,0},\OVV g^{\,1},\LD,\OVV g^{\,q-1}\}.
\ED

For an integer $k$ let $(k)_p$ be the representative of $k$ mod $p$, i.e., the number $j\in\{0,1,\LD, p-1\}$ fulfilling $j\equiv k$ mod p.
Let $b\ge 2$  be a natural number, $p\NDIV b$. We define the function
\BE
\label{2.1}
  \theta_b: \Z\to \Z: k\mapsto \theta_b(k)=\frac{b(k)_p-(bk)_p}{p}.
\EE
It is easy to illustrate the effect of this function. Indeed, for $k\in\{1,\LD,p-1\}$,
\BD
  \theta_3(k)=\begin{cases}
               0 & \hbox{ if $k<p/3$};\\
               1 & \hbox{ if $p/3\le k<2p/3$};\\
               2 & \hbox{ if $2p/3\le k$}.
               \end{cases}
\ED
We also need the property
\BE
\label{2.2}
  \theta_b(k)+\theta_b(p-k)=b-1 \MB{ for }k \in \{1,2,\LD,p-1\},
\EE
which is easy to check.
For a class $C\in G/H$ we define
\BD
  \SB_C=\sum_{k\in C}\theta_b(k).
\ED
In the case $b=2$, say, we count, thereby, the elements of $C$ that are $>p/2$.
Because of (\ref{2.2}) we have
\BD
  \sum_{C\in G/H}\SB_C=\sum_{1\le k\le p-1}\theta_b(k)=\sum_{1\le k<p/2}\theta_b(k)+\sum_{1\le k<p/2}\theta_b(p-k)=\frac{(b-1)(p-1)}2.
\ED
Accordingly, each sum $\SB_C$ has the {\em expected value}
\BE
\label{2.3}
  \EB=\frac{(b-1)(p-1)}{2q}.
\EE

The deviation of $\SB_C$, $C\in G/H$, from the expected value is measured by
\BD
  \TB_C=\SB_C-\EB,\enspace C\in G/H.
\ED
It turns out, in this context, that we can dispense with half of the classes $C$.
We write
\BD
  -C=\OVV{-1}C=\{p-k; k\in C\}
\ED
and obtain
\BE
\label{2.4}
  \TB_{-C}=-\TB_C,\enspace C\in G/H
\EE
from (\ref{2.2}).
Hence it suffices to  select one class from each pair $C, -C$.
The following choice is self-evident.
As above, let $g$ be a primitive root mod $p$.
Then the classes $\OVV g^{\,k}$ with even exponents $k\in\{0,\LD,q-1\}$ turn into those with odd exponents
by multiplication with $\OVV{-1}$.
In the sequel put $n=q/2$ and
\BD
  C_j=\OVV g^{2(j-1)}, j=1,\LD,n.
\ED
All these classes $C_j$ consist of quadratic residues mod $p$.
Observe that the order of the classes in $C_j$, $j\ge 2$, depends on the choice of the primitive root $g$ if $n\ge 3$.

Moreover, we write $\SB_j=\SB_{C_j}$ and $\TB_j=\TB_{C_j}$, $j=1,\LD, n$.
We call
\BD
  \TB=(\TB_1, \LD, \TB_n)
\ED
the {\em $b$-deviation vector} of $G/H$.

By the interrelation between the $b$-deviation vector and relative class numbers we will, in the main, obtain insight into two quantities connected with $\TB$:
on the one hand, the Euclidean norm $\|\TB\|$, i.e., the {\em dispersion} of the numbers $\SB_j$, $j=1,\LD,n$;
on the other hand, the {\em signs} of the numbers $\TB_j$, which show whether the sums $\SB_j$ are greater or smaller than their expected value $\EB$.
The assertions (\ref{1.1}), (\ref{1.2}), (\ref{1.3}), and (\ref{1.4}) emerge from these quantities in the case $q=2$.

\section*{3. Characters and Bernoulli numbers} 

We consider the $\C$-vector space $\C^n$, $n=q/2$, with the standard scalar product defined by
\BD
   \SP{z,u}=\sum_{j=1}^n z_j\OV{u_j},
\ED
where the bar denotes complex conjugation. A confusion with the bar used for residue classes is hardly possible.

We also consider the character group $X$ of $G/H$, whose elements
$\chi$ are understood as Dirichlet characters mod $p$. In particular, $\chi(k)=1$ if $k\in H$.
Since $\OVV{-1}\ne \OVV 1\in G/H$, one half of $X$ consists of odd characters, which are characterized by $\chi(-1)=-1$. Let $X^-$ denote the set of odd characters in $X$. We interpret each
$\chi\in X^-$ as a vector in $\C^n$ on identifying $\chi$ with $(\chi(C_1)), \LD,\chi(C_n))$.
The well-known orthogonality relations between characters show, for $\chi,\chi'\in X^-$,
\BD
  \SP{\chi,\chi'}=\begin{cases} 0 & \hbox{ if $\chi\ne\chi'$},\\
                                n & \hbox{ otherwise.}
                  \end{cases}
\ED
Hence the characters $\chi\in X^-$ form an orthogonal basis of $\C^n$. Therefore,
\BE
\label{3.2}
   z=\frac 1n\sum_{\chi\in X^-}\SP{z,\chi}\chi \:\MB { and }\: \|z\|^2=\frac 1n \sum_{\chi\in X^-}|\SP{z,\chi}|^2
\EE
for all $z\in \C^n$.

For a character $\chi\in X^-$,
\BD
  \BCH=\frac 1p \sum_{k=1}^{p-1}\chi(k)k
\ED
is the corresponding Bernoulli-number. Our first main result is as follows.

\begin{theorem} 
\label{s1}
For the $b$-deviation vector $\TB$ and each character $\chi\in X^-$,
\BD
  \SP{\TB,\chi}=(b-\chi(b))\BCHB/2,
\ED
where $\OV{\chi}$ denotes the complex conjugate character of $\chi$.
\end{theorem} 

\MN
{\em Proof.} By (\ref{2.4}), we have, since $\chi$ is odd,
\BD
  2\SP{\TB,\chi}=\sum_{C\in G/H}\TB_C\OV{\chi(C)}.
\ED
This equals
\BD
  \sum_{C\in G/H}\SB_C\OV{\chi}(C)-\EB\sum_{C\in G/H}\OV{\chi}(C).
\ED
Here the second sum vanishes since $\chi$ is a nontrivial character of $G/H$.
Therefore, (\ref{2.1}) gives
\BD
  2\SP{\TB,\chi}=\frac 1p \sum_{k=1}^{p-1}\theta_b(k)\OV{\chi}(k)=\frac bp\sum_{k=1}^{p-1}k\OV{\chi}(k)-\frac 1p\sum_{k=1}^{p-1}(bk)_p\OV{\chi}(k).
\ED
Only the last of these sums requires some discussion. Since $p\NDIV b$, the numbers $(bk)_p$ run through all $l\in\{1,\LD, p-1\}$. However, if $l=(bk)_p$, then $\OV k=\OV l\cdot\OV b\,^{-1}\in G$.
The assertion follows from $\OV{\chi}(\OV b\,^{-1})=\chi(b)$.
\STOP

\section*{4. Quadratic residues} 

It is easy to prove the identities (\ref{1.2}) and (\ref{1.3}) with the above tools. We will also prove another result of this kind; see (\ref{4.16}).

Let $q=2$, in particular, $p\equiv 3$ mod 4. Furthermore, let $p>3$.
Here $n=1$ and $C_1=C=H$, the set of quadratic residues mod $p$ of Section \ref{s1}. The only odd character of $G/H$ is the
Legendre symbol $\LS{-}{p}$, for which we write $\chi_2$ in the sequel. Moreover, $-B_{\chi_2}$ is the class number of $\Q(\sqrt{-p})$,
which is also the relative class number of this imaginary quadratic number field and denoted by
$h_2^-$; see \cite[p. 38]{Wa}.
The $b$-deviation vector $\TB$ has the form
\BD
 \TB=\TB_C=\SB_C-(b-1)(p-1)/4.
\ED
Because of (\ref{3.2}) and Theorem \ref{s1}, we have
\BE
\label{4.1}
  \TB_C=-(b-\chi_2(b))h_2^-/2
\EE
Accordingly,
\BE
\label{4.2}
   \SB_C= (b-1)(p-1)/4-(b-\chi_2(b))h_2^-/2.
\EE
In Section 2 we observed that, for $b=2$,
\BE
\label{4.3}
 \SZ_C=|C\cap[p/2,p]|.
\EE
Since $|C|=(p-1)/2$, we obtain
\BE
\label{4.4}
  |C\cap[0,p/2]|=(p-1)/2-\SZ_C.
\EE
Let $N=-C$ be as in Section 1. Then
\BD
\label{4.6}
 |N\cap[0,p/2]|=(p-1)/2-|C\cap[0,p/2]\}|=\SZ_C.
\ED
Together with (\ref{4.2}), this yields (\ref{1.2}), provided that we observe that $\LS 2p=1$ only if $p\equiv 7$ mod 8.

\MN
{\em Remark.}
From (\ref{4.2}) we see that the deviation  of $\SB_C$ from the expected value $\EB$ is small relative to $\EB=(b-1)(p-1)/4$ if $b$ is fixed
and $p$ tends to infinity. Indeed, the inequality of P\'olya-Vinogradov
implies $h_2^-\ll \sqrt p\log p$; see \cite[pp. 45, 214]{Wa}.

\MN
Next we introduce the function $f=1+\theta_2+\theta_3-\theta_6$.
It is easy to check that for $k\in\{1,\LD,p-1\}$
\BD
  f(k)=\begin{cases}
       1  & \hbox{ if $k<p/6$;}\\
       -1 & \hbox{ if $k>5p/6$;}\\
       0  & \hbox{ otherwise}.
       \end{cases}
\ED
We form
\BD
  S_C^f=\sum_{k\in C}f(k)=|C\cap[0,p/6]|-|C\cap[5p/6,p]|.
\ED
The expected value of this sum is
\BD
 E^f=|C|+E^{(2)}+E^{(3)}-E^{(6)}=0.
\ED
Hence we obtain
\BE
\label{4.13.1}
  S_C^f=S_C^f-E^f=T_C^{(2)}+T_C^{(3)}-T_C^{(6)}.
\EE
For $\chi_2=\LS{-}{p}$ we have, because of (\ref{4.1}),
\BE
\label{4.14}
 S_C^f=\SP{S_C^f,\chi_2}=\SP{T_C^{(2)}+T_C^{(3)}-T_C^{(6)},\chi_2}=-\WT c_{\chi_2}h_2^-/2,
\EE
with $\WT c_{\chi_2}=-1-\chi_2(2)-\chi_2(3)+\chi_2(6)$. But the number $\WT c_{\chi_2}$ takes the value $2$ if $\chi_2(2)=\chi_2(3)=-1$, and $-2$, otherwise.
Therefore, (\ref{4.13.1}) and (\ref{4.14}) yield the identity (\ref{1.3}).

One can also find $|C\cap[p/6,p/3]|-|C\cap[2p/3,5p/6]|$ if one uses the function $-\theta_2-2\theta_3+\theta_6$ instead of $f$, and
$|C\cap[p/3,p/2]|-|C\cap[p/2,2p/3]|$ by means of the function $-2\theta_2+\theta_3$; see \cite{Be}.

Finally, let $b=p+1$. The definition (\ref{2.1}) shows
$\theta_{p+1}(k)=k$ for $k\in\{1,2,\LD,p-1\}$ and
\BD
  \SPP_C=\sum_{k\in C} k.
\ED
Hence we call $\SPP_C$ the {\em class sum} of $C$.
By (\ref{4.2}), and since $\chi_2(p+1)=1$,
\BD
  \sum_{k\in C} k=p(p-1)/4-ph_2^-/2.
\ED
This shows that the class sum of $C$ is smaller than its expected value
$p(p-1)/4$.  As we remarked in the context of (\ref{1.2}), the elements of $C$ are accumulated in the lower half of the interval
$[0,p]$.
This convenes with the magnitude of the class sum. With $N=-C$ as in Section 1 we have
\BE
\label{4.16}
 \SPP_C-\SPP_N=-ph_2^-.
\EE

As concerns the distribution of quadratic residues by means of Dirichlet's formula
(\ref{1.3}), see also \cite{ChRo}. Further results on
the distribution of quadratic residues and higher power residues can be found in
\cite{Sh}, \cite{WaFa}, and \cite{WaXu}.

\section*{5. Digit expansions} 

Let $p\equiv 3$ mod 4, $p>3$ and $H=C_1=C$ be as in the foregoing section.

We will show that (\ref{1.4}) is a corollary to a more general theorem. To this end let
$b$ be a quadratic nonresidue mod $p$. Let $d$ be the order of the subgroup $\SP{\OV b}$ of $G$ (observe that $d$ is even). Moreover, let
$g$ be a primitive root  mod $p$ and $b_j=(g^{2j})_p$, $j=0,\LD,(p-1)/d-1$.  By \cite[Satz 1]{Gi3}, $b_j/p$ has the digit expansion
\BE
\label{4.12}
   b_j/p=\sum_{k=1}^{\infty}a_k^{(j)}b^{-k} \:\MB{ with }\: a_k^{(j)}=\theta_b(b_jb^{k-1}),
\EE
and the period $(a_1^{(j)}, a_2^{(j)}, \LD, a_d^{(j)} )$.

\begin{theorem} 
\label{s2}

In the above setting, let
\BD
  S_j=a_1^{(j)}+ a_3^{(j)}+ \LD+ a_{d-1}^{(j)},\enspace j=0,\LD,(p-1)/d-1,
\ED
be the sum of the digits with odd indices in the respective period.
Then
\BE
\label{4.12.2}
   S_0+S_1+\LD+S_{(p-1)/d-1}=(b-1)(p-1)/4-(b+1)h_2^-/2.
\EE

\MN
\end{theorem} 

\MN
{\em Proof.} The sum on the left-hand side of (\ref{4.12.2}) consists of the summands $\theta_b(b_jb^{k-1}_p)$, $j=0,\LD,(p-1)/d-1$, $k=1,3,\LD, d-1$.
But the residue classes $\OV b^{k-1}$, $k=1,3,\LD,d-1$ run through the group $H\cap\SP{\OV b}$ and the numbers $b_j$, $j=0,\LD,(p-1)/d-1$, are a system of representatives of
$H/H\cap\SP{\OV b}$.
Accordingly, the numbers $(b_jb^{k-1})_p$ run through all elements of $H=C$. So this sum equals $\SB_C$.
Since $\chi_2(b)=-1$, the assertion follows from  (\ref{4.2}).
\\ \MB{ }\STOP

\MN
If $b$ is a primitive root mod $p$, i.e., if $d=p-1$, only $b_0=1$ occurs in Theorem \ref{s2}. Then $b_0/p=1/p$ has the period $(a_1,a_2,\LD,a_{p-1})$. The theorem
says
\BE
\label{4.8}
 a_1+a_3+\LD+a_{p-2}=(b-1)(p-1)/4-(b+1)h_2^-/2.
\EE
From (\ref{2.2}) we derive
\BD
  a_{k}+a_{k+(p-1)/2}=b-1, \enspace k=1,\LD,p-1.
\ED
But  $a_{k+(p-1)/2}$ runs through the digits $a_2, a_4,\LD, a_{p-1}$, if $k$ runs through $1,3,\LD, p-2$.
Hence we have
\BD
  (a_1+a_3+\LD+a_{p-2})+(a_2+a_4+\LD+a_{p-1})=(b-1)(p-1)/2.
\ED
Together with (\ref{4.8}), this yields (\ref{1.4}).

\MN
{\em Example.} A simple example of Theorem \ref{s2} is $p=11$ and $b=10$; so $d=2$. With $2$ as a primitive root mod $p$, we obtain
$b_0=1$, $b_1=4$, $b_2=5$, $b_3=9$, $b_4=3$. We have $1/p=0.\OV{09}$, $4/p=0.\OV{36}$, $5/p=0.\OV{45}$, $9/p=0.\OV{81}$, and $3/p=0.\OV{27}$.
Hence the sum on the left-hand side of (\ref{4.12.2}) is $0+3+4+8+2=17$. On the other hand, the right hand side of (\ref{4.12.2}) is $9\cdot10/4-11\cdot h_2^-/2=17$ since $h_2^-=1$.

\MN
Next let $b$ be a quadratic residue mod $p$. So $\SP{\OV b}$ is a subgroup of $H$. Let $|\SP{\OV b}|=d$.
As above, let $g$ be a primitive root  mod $p$ and define $b_j=(g^{2j})_p$, $j=0,\LD,(p-1)/(2d)-1$. We have
\BD
   b_j/p=\sum_{k=1}^{\infty}a_k^{(j)}b^{-k} \:\MB{ with }\: a_k^{(j)}=\theta_b(b_jb^{k-1}),
\ED
and the period $(a_1^{(j)}, a_2^{(j)}, \LD, a_d^{(j)} )$.

\begin{theorem} 
\label{s2.1}

In the above setting, let
\BD
  S_j=a_1^{(j)}+ a_2^{(j)}+ \LD+ a_d^{(j)}
\ED
be the sum of the digits of the period of $b_j/p$, $j=0,\LD,(p-1)/(2d)-1$. Then
\BE
\label{4.12.1}
   S_0+S_1+\LD+S_{(p-1)/(2d)-1}=(b-1)(p-1)/4-(b-1)h_2^-/2.
\EE

\MN
\end{theorem} 

\MN
The {\em proof} is quite analogous to the proof of Theorem \ref{s2} and, therefore, omitted.

In the case $\SP{\OV b}=H$, i.e., $d=(p-1)/2$ we have only the digit sum $S_0$ of the period of $1/p$ on the left-hand side of (\ref{4.12.1}). This special case of
Theorem \ref{s2.1} was shown in \cite{Gi3}; see also \cite{MuTh}.

\MN
{\em Example.} Let $p=79$ and $b=10$, a quadratic residue mod $p$. We have $|\SP{\OV{10}}|=13$ and $(p-1)/(2d)=3$.
We take $g=3$ and obtain $b_0=1$, $b_1=9$, $b_2=2$.
Now
\BD
 1/79=0.\OV{0126582278481},\: 9/79=0.\OV{1139240506329},\: 2/79=0.\OV{0253164556962},
\ED
and the digit sums of the respective periods are $54$, $45$, $54$. Their sum is $153$.
On the other hand, the right-hand side of (\ref{4.12.1}) is $9\cdot 78/4-9\cdot h_2^-/2=153$ since $h_2^-=5$.

\begin{theorem} 
\label{s2.2}
Let $p\equiv 3$ mod 4, $p>3$, and
$b\ge 2$, $p\NDIV b$. In addition, let $b$ be even and of the order $(p-1)/2$ mod $p$. Then we have, for the period $(a_1,\LD,a_{(p-1)/2})$ of $1/p$ ,
\BE
\label{4.13}
  |\{k; a_k\le b/2-1\}|-|\{k; a_k\ge b/2\}|=(2-\chi_2(2))h_2^-.
\EE

\end{theorem} 

\MN
{\em Remark.} Whereas $b$ appears on the left-hand side of (\ref{4.13}), it does not appear on the right hand side.
Accordingly, the difference on the left-hand side is the same for all $b$ in question.

\MN
{\em Proof of Theorem} \ref{s2.2}. The digits of the period have the form $\theta_b(j)$, $j\in C$.
One shows, for each $j\in C$, the equivalence of two assertions, namely,
\BD
   \theta_b(j)\le b/2-1 \:\MB{ and }\: j\le p/2.
\ED
For this purpose one uses the estimates
\BD
 \frac{bj-1}p\ge \frac{bj-(bj)_p}p\ge \frac{bj+1-p}p.
\ED
Then the theorem follows from (\ref{4.2}), (\ref{4.3}) and (\ref{4.4}).
\STOP

\MN
{\em Remark.} The adaption of Theorem \ref{s2.2} to the setting of Theorem \ref{s2.1} is left to the reader.

\MN
The connection between digit expansions and class numbers is also investigated in \cite{ChKr}, \cite{Hi}, \cite{Mi}, and \cite{PuSa}.
As an example of an analytical result on the digit expansion of $m/p$, $(m,p)=1$, we mention \cite[Cor. 8]{BoKoSh}.

\section*{6. Sixth power residues} 

In the case $q=2$, $n=1$, the $b$-deviation vector $\TB=\TB_C$ satisfies the equation (\ref{4.1}), whose right-hand side is always negative since $\chi_2(b)\in\{\pm 1\}$. Accordingly, the
sum $\SB_C$ is always smaller than its expected value. In the context of (\ref{1.4}), the sum of the digits with odd indices is smaller than the sum of the digits with even indices
--- which is a consequence of the same fact.
Therefore, we think that it is desirable  to know the signs of the numbers  $\TB_j$, $j=1,\LD,n$, for greater values of $q=2n$.
We will see that, under certain conditions, this is possible for $q=6$.

For the time being, we assume $q=2n$, where $n$ is odd.
Recall the setting of Sections 2, 3. In particular,
$p\equiv q+1$ mod $2q$, $p>q+1$. For a primitive root $g$ mod $p$ we have $C_j=\OVV g^{2(j-1)}$, $j=1,\LD, n$. For $b\ge 2$,
we consider $\SB_j=\SB_{C_j}$ and $\TB_j=\TB_{C_j}$, $j=1,\LD, n$.
The $b$-deviation vector is $\TB=(\TB_1,\LD,\TB_n)$.

The set $X^-$ of odd characters of $G/H$ satisfies $|X^-|=n$.
Since $n$ is odd, we have $p\equiv 3$ mod 4. Therefore, $\chi_2 =\LS{-}p\in X^-.$  In the sequel we use the abbreviation
\BD
   c_{\chi}=b-\chi(b), \enspace \chi\in X^-.
\ED
Then the assertion of Theorem \ref{s1} takes the form
\BE
\label{5.2}
  \SP{\TB,\chi}=c_{\chi}B_{\OV{\chi}}/2.
\EE

Because of (\ref{4.1}) and $\chi_2(C_{j})=1$, $j=1,\LD,n$, the vector $\TB$ lies
in the hyperplane
\BE
\label{5.3}
  P=\{(z_1,\LD,z_n)\in\R^n; z_1+\LD+z_n=-c_{\chi_2}h_2^-/2\}\enspace (\subseteq \R^n).
\EE

Let the sign of a real number $x$ be defined in the usual way, namely,
\BD
  \SIGN(x)=\begin{cases}
            1& \hbox{ if } x>0;\\
            0& \hbox{ if } x=0;\\
           -1& \hbox{ if } x<0.
           \end{cases}
\ED
For $z=(z_1,\LD,z_n)\in \R^n$ we define the {\em sign vector}
\BD
  \SIGN(z)=(\SIGN(z_1),\LD,\SIGN(z_n)).
\ED
The point of $P$ with the smallest (Euclidean) norm has the form
\BD
  z^{(0)}=-c_{\chi_2}h_2^-(1,1,\LD,1)/n.
\ED
Its sign vector is $(-1,\LD,-1)$, which we call the {\em main type}. The following idea is plausible: Namely,
the closer $\TB$ is to $z^{(0)}$, the more $\SIGN(\TB)$ resembles the main type.
We have the following theorem.

\begin{theorem} 
\label{s3}
Let $n>1$ be odd.
\\ {\rm (a)} The vector $\SIGN(\TB)$ agrees with the main type in at least one position.
\\ {\rm (b)} If
\BD
 \|\TB \|^2\ge c_{\chi_2}^2\frac{(h_2^-)^2}4,
\ED
then $\SIGN(\TB)$ differs from the main type.
\\ {\rm (c)} Let $k\in\{1,2,\LD,n-1\}$. If
\BD
  \|\TB \|^2< c_{\chi_2}^2\frac{(h_2^-)^2}{4k},
\ED
then $\SIGN(\TB)$ agrees with the main type in at least $k+1$ positions.

\end{theorem} 

\MN
The theorem is an immediate consequence of the following.

\begin{lemma} 
\label{l1}
Let $n>1$ and
$x_1,\LD,x_n\in \R$ such that $x_1+\LD+x_n=1$.
\\ {\rm (a)} There is a $j\in\{1,\LD,n\}$ such that $x_j>0$.
\\ {\rm (b)} If $x_1^2+\LD+x_n^2\ge 1$, then there is a $j$ with $x_j\le 0$.
\\ {\rm (c)} Let $k\in\{1,\LD,n-1\}$. If $x_1^2+\LD+x_n^2<1/k$, there are at least $k+1$ indices $j$ with $x_j>0$.

\end{lemma} 

\MN
{\em Proof.} Assertion (a) is obvious.
\\ Ad (b): If all $x_j>0$, then $n> 1$ implies $0<x_j<1$ for all $j$. Thus,
$x_1^2+\LD+x_n^2<x_1+\LD+x_n=1$.
\\ Ad (c): We write $x=(x_1,\LD,x_n)$ and use the standard scalar product and the corresponding Euclidean norm in $\R^n$.
Suppose $|\{j; x_j>0\}|\le k$. Then we may assume, without loss of generality, $x_{k+1},\LD, x_n\le 0$.
Let $(1,1,\LD,1,0,0,\LD,0)\in \R^n$ be the vector with exactly $k$ entries equal to $1$. We have
$1\le x_1+\LD+x_k=\SP{x,(1,1,\LD,1,0,0,\LD,0)}$. Cauchy's inequality shows $1\le\|x\|\sqrt k$.
\STOP

\MN
{\em Remark.} The Lemma is sharp, inasmuch one does not receive stronger statements from its premises.
In the case $n=3$, which is the most relevant for us, we will get a graphic demonstration that no stronger result is possible.

\MN
Now to the case $q=6$, i.e., $n=3$. Here $X^-=\{\chi_2,\chi_6,\OV{\chi_6}\}$,
where $\chi_6$ is a character of order 6 and, as above, $\OV{\chi_6}=\chi_6^{-1}$ its complex conjugate character.
Because of (\ref{3.2}) and (\ref{5.2}),
\BE
\label{5.4}
  \|\TB\|^2=\frac 13(|\SP{\TB,\chi_2}|^2+2|\SP{\TB,\chi_6}|^2)=\frac 13(c_{\chi_2}^2\frac{B_{\chi_2}^2}4+ 2|c_{\chi_6}|^2\frac{|B_{\chi_6}|^2}4).
\EE
The Bernoulli numbers in this formula are connected with relative class numbers. To this end let $\Q^{(p)}$ be the $p$th cyclotomic field. We identify its Galois group over $\Q$ in
the usual way with the cyclic group $G=(\Z/p\Z)^{\times}$.
Since $p\equiv 7$ mod 12, the field $\Q^{(p)}$ has exactly one subfield $K_6$ of degree $[K_6:\Q]=6$. The field $K_6$ is the fixed field of
the subgroup $H\subseteq G$ of index 6; because $\OV{-1}\not\in H$, $K_6$ is imaginary.
Let $h_6^-$ be its relative class number.  By
\cite[Th. 4.17]{Wa} and \cite[Satz 23]{Ha}, we have $h_6^-=|B_{\chi_2}B_{\chi_6}^2|/4$ (observe $p>7$). Moreover,  $K_6$ contains a unique subfield $K_2$ of degree $[K_2:\Q]=2$, the fixed field of
$H'=\{\OV k^2; \OV k\in G\}$.
Because $p\equiv 3$ mod 4, $\OV{-1}\not\in H'$
and, thus, $K_2$ is also imaginary. In addition, $K_2\subseteq K_6$.
The relative class number $h_2^-$ of  $K_2$ is given by $|B_{\chi_2}|=h_2^-$. With this notation, (\ref{5.4}) yields the following.

\begin{theorem} 
\label{s4}
Let $p\equiv 7$  mod $12$, $p>7$. Let $q=6$ and $b\ge 2$, $p\NDIV b$. The $b$-deviation vector $\TB$ satisfies
\BE
\label{5.6}
   \|\TB\|^2=\frac 13(c_{\chi_2}^2\frac{(h_2^-)^2}4+ 2|c_{\chi_6}|^2\frac{h_6^-}{h_2^-}).
\EE
\end{theorem} 

\MN
So the theorem expresses the {\em variance} of the sums $S_j$, $j=1,2,3$, by relative class numbers.
Combined with Theorem \ref{s3} and some computation, formula (\ref{5.6}) gives the following corollary.

\begin{corollary} 
\label{k1}
In the setting of Theorem {\rm \ref{s4}}, the following assertions hold.
\\ {\rm (a)} The sign vector $\SIGN(\TB)$ agrees with the main type $(-1,-1,-1)$ in at least one position.
\\ {\rm (b)} If
\BD
   |c_{\chi_6}|^2h_6^-\ge c_{\chi_2}^2\frac{(h_2^-)^3}4,
\ED
then $\SIGN(\TB)$ differs from the main type.
\\ {\rm (c1)} If
\BD
   |c_{\chi_6}|^2h_6^-< c_{\chi_2}^2\frac{(h_2^-)^3}4,
\ED
then $\SIGN(\TB)$ agrees with the main type in at least two positions.
\\ {\rm (c2)} If
 \BD
   |c_{\chi_6}|^2h_6^-< c_{\chi_2}^2\frac{(h_2^-)^3}{16},
\ED
then $\SIGN(\TB)$ is the main type.

\end{corollary} 

\MN
In the case $q=6$, the group $H$ consists of the sixth powers in $G$ and equals the class $C_1$. Moreover, $C_2$ and $C_3$ are the cosets of $H$ in the group $H'$ of the squares in $G$.
These cosets have the form
$\OVV{g}^{\,2}$ and $\OVV{g}^{\,4}$, respectively; their order depends on the choice of the primitive root $g$ mod $p$.

\MN
{\em Example.}
We compare the cases $b=2$ and $b=p+1$. In the case $b=2$,
\BD
  \TB_j=|C_j\cap[p/2,p]|- (p-1)/12,
\ED
in the case $b=p+1$, however,
\BD
 \TB_j=\sum_{k\in C_j} k-p(p-1)/12,
\ED
$j=1,2,3$; see (\ref{2.3}).

For the primes $p\equiv 7$ mod 12, $7<p<500$, the situation is as follows. For 13 of 23 primes in question the criterion (b) of Corollary \ref{k1} applies both for
$b=2$ and $b=p+1$. For the remaining ten primes the respective information is given in Table 1.

\MN
\begin{center}
\begin{tabular}{c|c|c|c|c|c|c|c|c|c} 

 $p$ & $g$ & $c_6$ & $c_2$ & $h_6^-$ & $h_2^-$ &  $b=2$ & $\SIGN$    & $b=p+1$ & $\SIGN$    \\ \hline\rule{0mm}{5mm}
 79  &  3  &     7 &     1 &      5  &      5  &  (b)   & $-1,-1,1$  &  (c2)   & main type \\ \rule{0mm}{5mm}
 103 &  5  &     7 &     1 &      5  &      5  &  (b)   & $1,-1,-1$  &  (c2)   & main type    \\ \rule{0mm}{5mm}
 139 &  2  &     3 &     9 &      9  &      3  &  (c1)  & main type  &  (b)    & $1,-1,-1$  \\  \rule{0mm}{5mm}
 151 &  6  &     7 &     1 &      7  &      7  &  (c1)  & $1,-1,-1$  &  (c2)   & main type   \\ \rule{0mm}{5mm}
 199 &  3  &     7 &     1 &     27  &      9  &  (b)   & $1,-1,-1$  &  (c2)   & main type   \\ \rule{0mm}{5mm}
 271 &  6  &     7 &     1 &     11  &     11  &  (c2)  & main type  &  (c2)   & main type   \\ \rule{0mm}{5mm}
 367 &  6  &     7 &     1 &     27  &      9  &  (b)   & $-1,-1,1$  &  (c2)   & main type   \\ \rule{0mm}{5mm}
 439 & 15  &     1 &     1 &    405  &     15  &  (c1)  & $-1,-1,1$  &  (c1)   & $-1,-1,1$  \\  \rule{0mm}{5mm}
 463 & 3   &     7 &     1 &     49  &      7  &  (b)   & $1,-1,-1$  &  (c1)   & $1,-1,-1$  \\  \rule{0mm}{5mm}
 487 & 3   &     7 &     1 &     49  &      7  &  (b)   & $-1,-1,1$  &  (c1)   & $-1,-1,1$

\end{tabular} 

\MN
{\bf Table 1}
\end{center}

\MN
In the first line of the table we find $c_6=|c_{\chi_6}|^2$ and $c_2=c_{\chi_2}^2$ for the case $b=2$. In the case $b=p+1$, we have $|c_{\chi_6}|^2=p^2=c_{\chi_2}^2$, hence these constants play no role in the
application of the criteria of Corollary \ref{k1}. In the columns headed $\BQ b=2\EQ$ and $\BQ b=p+1\EQ$ one finds the criterion that may be applied;
here \BQ(c1)\EQ says that the criterion (c2) cannot be applied.
In the columns labeled \BQ $\SIGN$\EQ we omit the brackets of the sign vectors; for instance,  $-1,-1,1$ stands for $(-1,-1,1)$. In one case, namely $p=139$, $b=2$,
we have the main type although only (c1) is applicable.

Among the 13 primes $p\equiv 7$ mod 12, $7<p<500$, not rendered in the table we find sign vectors with two entries 1, such as
$p=163$, $g=2$, with $(1,1,-1)$ in the case $b=p+1$.

Observe that Corollary \ref{k1} can be applied to
\BD
  |C_j\cap[0,p/6]|-|C_j\cap[5p/6,p]=T_j^{(2)}+T_j^{(3)}-T_j^{(6)},\enspace j=1,2,3,
\ED
where the constant $c_{\chi_6}$ has to be replaced by $\WT c_{\chi_6}= -1-\chi_6(2)-\chi_6(3)+\chi_6(6)$ and
$c_{\chi_2}$ by $\WT c_{\chi_2}= -1-\chi_2(2)-\chi_2(3)+\chi_2(6)$. Since $\WT c_{\chi_2}=\pm 2$, the main type takes two possible forms, namely $(-1,-1,-1)$ and $(1,1,1)$, in the respective cases; see (\ref{4.14}).

Corollary \ref{k1} can also be applied to {\em digit expansions}. As an example, we assume that $b$ has the order $(p-1)/2$ mod $p$, a case considered in connection with Theorem \ref{s2.1}. Then $1/p$ has the period $(a_1,\LD,a_{(p-1)/2})$. As above, we have three classes $C_{j}=\OVV b^{\,j-1}$, $j=1,2,3$. In this situation, the sums
$\SB_{j}$ satisfy
\BD
  \SB_{j}=\sum_{k=0}^{(p-1)/6-1}a_{j+3k}\: \MB{ and }\: \TB_j=\SB_{j}-(b-1)(p-1)/12,\enspace j=1,2,3.
\ED
In the case $b=10$, we have $|c_{\chi_6}|^2=111$ and $c_{\chi_2}^2=81$. There are 15  primes $p\equiv 7$ mod 12, $7<p<1000$, such that 10 has the order $(p-1)/2$ mod $p$. For 10 of these primes the criterion
(b) of Corollary \ref{k1} applies. Table 2 describes the situation for the remaining five primes. It is organized like Table 1.

\MN
\begin{center}
\begin{tabular}{c|c|c|c|c} 

$p$ & $h_6^-$  & $h_2^-$ &  $b=10$ & $\SIGN$ \\ \hline\rule{0mm}{5mm}

 151 &       7 &       7 &  (c2)   & main type \\ \rule{0mm}{5mm}
 199 &      27 &       9 &  (c2)   & main type \\ \rule{0mm}{5mm}
 439 &     405 &      15 &  (c1)   & $-1,1,-1$ \\ \rule{0mm}{5mm}
 631 &     325 &      13 &  (c1)   & $-1,1,-1$  \\ \rule{0mm}{5mm}
 991 &     119 &      17 &  (c2)   & main type

\end{tabular} 

\MN
{\bf Table 2}
\end{center}

\MN
Next we elucidate Lemma 1 in the case $n=3$. We write $x=(x_1,x_2,x_3)$ for a vector in $\R^3$. Then
\BD
  P=\{x\in\R^3;x_1+x_2+x_3=1\}
\ED
is the plane spanned by the standard base vectors $e_1=(1,0,0)$, $e_2=(0,1,0)$ and $e_3=(0,0,1)$. This plane is rendered in Diagram 1.

Here the points $e_j$ are the corners of the equilateral triangle in the middle, whose sides have the length $\sqrt 2$.
The three straight lines through these points divide
$P$ into seven zones. The diagram shows the value of $\SIGN(x)$ for each zone.
Whereas we have $(1,1,1)$ (as the main type) for the points in the interior of the triangle, we have one or two negative signs in  the interior of the other zones.
The circle through $e_1, e_2, e_3$ consists of the points
$x\in P$ with $\|x\|^2=1$. It is the circumcircle of the triangle. On this circle and outside it we have the situation of criterion (b) of the lemma. Indeed, the diagram shows that the main type does not occur
there. Criterion (c) of the lemma with $k=1$ (which corresponds to criterion (c1) of
Corollary \ref{k1}) applies to the points $x$ in the interior of this circle. For these points $\SIGN(x)$ has at least two entries 1.

The circle through the points $(e_1+e_2)/2$, $(e_1+e_3)/2$, $(e_2+e_3)/2$ consists of all points $x\in P$ with $\|x\|^2=1/2$ and is the incircle of the triangle.
In the interior of the incircle we have $\SIGN(x)=(1,1,1)$,
which corresponds to criterion (c2) of Corollary \ref{k1}.
The diagram also shows that a better criterion for $\SIGN(x)$ based only on
$\|x\|$ does not exist, since all other circles concentric with these two circles are inferior.
In particular, no such criterion for the occurrence of two signs $-1$ in $\SIGN(x)$ is possible.

\input{Power_Pict_1}

In view of the case $n=3$, it is not difficult to visualize the case $n=4$. Here the equilateral triangle is replaced by a regular tetrahedron.
The planes containing its sides give rise to 15 zones in the space $P$, where each zone is connected with a certain sign pattern. Instead of the two circles we have
three spheres, namely, through the corners, the midpoints of the edges, and the midpoints of the sides of the tetrahedron.

Suppose that $\SIGN(\TB)$  is the main type in the case $n=3$. Then criterion (c1) of Corollary \ref{k1} applies and says that $\SIGN(T)$ differs from the main type in at most one position.
The proportion of the area of the incircle to the area of the triangle in Diagram 1 suggests a chance of about 60 \% for the applicability of criterion (c2).
However, this is only true if the main type is uniformly distributed in the triangle. But it seems that this does not hold, see Diagram 2 and its explanations.

The $(p+1)$-deviation vector $T^{(p+1)}$ describes the difference between the {\em class sums}
\BD
  S^{(p+1)}_j=\sum_{k\in C_j}k, \enspace j=1,2,3,
\ED
and their expected value $p(p-1)/12$.
Diagram 2 shows the distribution of 4754 vectors $T^{(p+1)}$ for $7<p<3\cdot 10^5$, where $g$ is the smallest primitive root mod $p$.
These vectors have been {\em normalized}, i.e., divided by $-ph_2^-/2$, with the effect that the sum of the three components equals 1; see (\ref{5.3}).
For this reason we can adopt the scheme of Diagram 1 for our visualization. The main type corresponds to the sign vector $(1,1,1)$.
Only those vectors (1748 in number) whose image points fall outside our display window have been omitted.
The main type seems to occur with a frequency of about 25 \%  of all primes $p\equiv 7$ mod 12.  Computations suggest a chance of about 70 \% for the applicability of the criterion (c2) provided that the sign vector
is the main type. This relatively high percentage comes from the higher density of image points in the interior of the incircle of the triangle. But there are no image points near the midpoint of this circle.

This empty space can be explained as follows. Due to Theorem \ref{s4}, the normalized $(p+1)$-deviation vector has the number $(1+8h_6^-/(h_2^-)^3)/3$ as the square of its norm.
But the minimum of the numbers $h_6^-/(h_2^-)^3$ for primes in the range in question is about $0.0039$. This means that the said vector keeps some distance to the midpoint $(1,1,1)/3$ of the incircle.
Of course, this distance could become smaller if larger primes $p$ are included.

\input{Power_Pict_2}

The Brauer-Siegel theorem implies
\BD
  \log(h_6^-)/\log(( h_2^-)^3) \to 1
\ED
for $p\to\infty$; see \cite[pp. 42--45]{Wa}. This says that $h_6^-$ and $(h_2^-)^3$ have the same order of magnitude, but only very roughly.

The application of the criteria of Corollary \ref{k1} in the case just considered requires only the relative class numbers
$h_6^-$ and $h_2^-$ (since $c_{\chi_2}^2=|c_{\chi_6}^2|=p^2$). There are effective methods for the computation of these numbers; see \cite{Lo}.
Thereby, one could  find larger primes (say $p>10^{12}$) belonging to the main type.

It should be mentioned that for other values of $b$, in particular, for $b=2$, one obtains pictures of normalized $b$-deviation vectors that are fairly different from Diagram 2.

For sixth power residues see also \cite{SuZh}.

\section*{7. Tenth power residues} 

In the case $q=10$, let $p\equiv 11$ mod 20 and $p>11$. We restrict ourselves to $b=p+1$, i.e., $\SB_j=\sum_{k\in C_j}k$ is the class sum of $C_j=\OVV{g}^{\,2(j-1)}$, $j=1,\LD,5$, and $\EB$ equals $p(p-1)/20$.
We have $X^-=\{\chi_2\}\cup\{\chi_{10},\OV{\chi_{10}},\chi_{10}^3,\OV{\chi_{10}}^3\}$, where $\chi_{2}=\LS{-}p$ and $\chi_{10}$ is a character of order 10.
The analogue of (\ref{5.4}) has, since $b=p+1$, the form
\BE
\label{6.2}
  \|\TB\|^2=\frac{p^2}5({B_{\chi_2}^2}/4+ {|B_{\chi_{10}}|^2}/2+{|B_{\chi_{10}^3}|^2}/2).
\EE
The main difference between this case and the case $q=6$ consists in the fact that $|B_{\chi_{10}}|^2$ and $|B_{\chi_{10}^3}|^2$ cannot be expressed by rational class number factors. Indeed, we only have
\BD
   |B_{\chi_{10}}|^2|B_{\chi_{10}^3}|^2=16\, h_{10}^-/h_2^-,
\ED
where $h_{10}^-$ and $h_2^-$ are the relative class numbers of the subfields of $\Q^{(p)}$ of the degrees 10 and 2, respectively.
By the arithmetic-geometric inequality,
\BE
\label{6.4}
   \|\TB\|^2\ge \frac{p^2}5((h_2^-)^2/4+4\sqrt{h_{10}^-/h_2^-}).
\EE
Therefore, criterion (b) of Corollary \ref{k1} can be applied in the strict sense, however, no other one.
Indeed, if
\BD
  h_{10}^-\ge (h_2^-)^5/16,
\ED
then $\SIGN(\TB)$ differs from the main type.
Nevertheless, computations for $11<p<2\cdot 10^6$ show that the deviation of $\|\TB\|^2$ from the right-hand side of (\ref{6.4}) is small in most cases. This means that we have {\em probabilistic} criteria
like the following.
If
\BD
  h_{10}^-< (h_2^-)^5/4096,
\ED
then $\SIGN(\TB)$ is (probably) the main type. This procedure gives only two false predictions for $p<10^5$; see \cite[Sect. 3.3]{Gi2}, where further details can be found.

\medskip
\centerline{{\bf Funding}}

This research did not receive any specific grant from funding agencies in the public, commercial, or not-for-profit sectors.


\MN
Kurt Girstmair\\
Institut f\"ur Mathematik \\
Universit\"at Innsbruck   \\
Technikerstr. 13/7        \\
A-6020 Innsbruck, Austria \\
Kurt.Girstmair@uibk.ac.at

\end{document}